\documentclass[twoside,11pt]{article}
\usepackage{graphicx}
\usepackage{amssymb}
\usepackage{amsthm}
\newfont{\bbb}{msbm10 scaled\magstep 1}

\font\bigbold=cmbx10 at 14 pt
\font\bbigbold=cmbx10 at 17 pt 
\date{}

\title 
{\bbigbold An estimate of the centroid Banach-Mazur distance between planar convex bodies}

\begin{document}

\baselineskip 17pt 

\maketitle

\vskip -1.3cm
\centerline
{\bigbold Marek Lassak}

\vskip 0.15cm 
\pagestyle{myheadings} \markboth{Marek Lassak}{An estimate of the centroid Banach-Mazur distance}

\vskip 0.35cm
\noindent
{\bf Abstract.} 
We consider the variant of the Banach-Mazur distance $\delta_{BM}^{\rm cen} (C, D)$ of two convex bodies $C, D$ of $E^d$ with the additional requirement that the centroids of them coincide.
 We prove that $\delta_{BM}^{\rm cen} (C, D) \leq \frac{69}{17}$ for every convex bodies $C, D$ of $E^2$.



\date{}

\maketitle

\section{Introduction}

By ${\mathcal C}^d$ we denote the family of convex bodies (i.e., closed bounded convex sets with non-empty interior) of the $d$-dimensional Euclidean space $E^d$:

Our Theorem concerns the two-dimensional case of the so called {\it centroid Banach-Mazur distance} of convex bodies $C, D \in {\mathcal C}^d$:

$$\delta_{\rm BM}^{\rm cen} (C,D) = \inf_{a, h_\lambda} \{\lambda ; \, a(C) \subset D \subset h_\lambda a(C) \ {\rm and} \ {\rm cen}(a(C)) =  {\rm cen}(D) \},$$

\noindent
where $a$ stands for an affine transformation and $h_\lambda$ denotes a homothety with a ratio $\lambda \geq 1$ whose center is at the centroids ${\rm cen}(a(C)) =  {\rm cen}(D)$ of $a(C)$ and $D$. 
For the notion and properties of the centroid see \S2 and \S7 of the book \cite{[B+F]} by Bonnesen and Fenchel. 
For the notions of affine transformation and homothety see Part 1.1 of the book \cite {[T]} by Toth.

For every $C$ and $D$ we have $\delta_{\rm BM}^{\rm cen} (C,D) = \delta_{\rm BM}^{\rm cen} (D,C)$ (see Section 3).

In connection with this, recall that the {\it extended Banach-Mazur distance} of $C, D \in {\mathcal C}^d$ is defined as 

$$\delta_{\rm BM}(C, D) = \inf_{a,\; h_{\lambda}} \{ \lambda;\ \ a(D) \subset C \subset h_{\lambda}\big( a(D) \big) \},$$

\noindent
where $a$ stands for an affine transformation and $h_\lambda$ denotes a homothety with ratio $\lambda \geq 1$.
It is well known that this is a generalization of the original definition for centrally-symmetric $C$ and $D$ as given by Banach \cite{[Ba]} in behalf of him and Mazur (formally it is shown in Claim of \cite{[L3]}). 
The centroids of $C, D \in \mathcal{C}^d$ in the definition of $\delta_{\rm BM}^{\rm cen} (C,D)$ take over the roles of the centers of the centrally-symmetric bodies in the original definition of the Banach-Mazur distance.
A survey on the Banach-Mazur distance is given in the book \cite{[T-J]} by Tomczak-Jaegerman.
Moreover, in Sections 3.2 and 3.3 of the book \cite{[T]} by Toth, and in Section 4.1 of the book \cite{[AS]} by Aubrun and Szarek. 

As usual, by an {\it affine-regular hexagon} we understand a non-degenerated affine image of the regular hexagon.
Besicovitch \cite{[Be]} proved the following theorem

\vskip0.1cm
($\heartsuit$) {\it For every $A \in \mathcal{C}^2$ there is an affine-regular hexagon  $H_A$  inscribed in $A$.}

\vskip0.1cm
In our Theorem we prove that for every planar convex bodies $C$ and $D$ we have $\delta_{BM}^{\rm cen} (C, D) \leq \frac{69}{17}$, which is $4.058\dots$.
In the proof we substantially apply the following fact proved in \cite{[L4]}. 

\vskip0.1cm
($\diamondsuit$) {\it Let $A \in \mathcal{C}^2$ and $H_A$ be an affine-regular hexagon inscribed in $A$. 
Then the centroid of $A$ belongs to the homothetic image of $H_A$ with ratio $\frac{4}{21}$ and center in the center of $H_A$.}

\section{An estimate of the centroid BM-distance}

{\bf Theorem.}
{\it The centroid Banach-Mazur distance between arbitrary two planar convex bodies is at most $\frac{69}{17}$.}

\begin{proof}
Consider any convex bodies $C, D$ of $E^2$.
By ($\heartsuit$) there exists an affine-regular hexagon $H_C = c_1 \dots c_6$ inscribed in $C$ and an affine-regular hexagon $H_D = d_1 \dots d_6$ inscribed in $D$. 
Since we deal with $C$ and $D$, and thus with $H_C$ and $H_D$, up to non-degenerated affine images, we may assume that $c_i = d_i = (\cos \frac{2\pi}{6}i, \sin \frac{2\pi}{6}i)$, where $i=1, \dots, 6$,
in a rectangular coordinate system.

By ($\diamondsuit$) the centroid $(p,q)$ of $C$ and the centroid $(r^*,s^*)$ of $D$ are in $\frac{4}{21} H_C = \frac{4}{21} H_D$. 
Since this hexagon is the union of $12$ congruent triangles (obtained by dissecting the hexagon by its $12$ axes of symmetry), without loosing the generality we may assume that the centroids $(p,q)$ and $(r^*,s^*)$ are in the triangle $T$ with vertices 
$(0,0)$,  $(\frac{4}{21}, 0)$ and $(\frac{1}{7},\frac{1}{21}\sqrt 3)$.

Our considerations are not narrower if we assume that for instance $p \leq r^*$ (in the opposite case we change the roles of $C$ and $D$).

We provide the affine transformation $\tau$ given by $r= \frac{3}{2}r^* + \frac{1}{2}\sqrt 3 s^*$, $s= -\frac{1}{2}\sqrt 3 r^* + \frac{3}{2}s^*$ (it is the superposition of the rotation by $-30^\circ$ with a homothety with ratio $\sqrt 3$).
The obtained points $(r,s)$ form the triangle $T^+$ with vertices $(0,0)$, $(\frac{2}{7}, -\frac{2}{21}\sqrt 3)$ and $(\frac{2}{7}, 0)$.

Since $r^* \geq 0$ and $s^*\geq 0$, we have $r = \frac{3}{2}r^* + \frac{1}{2}\sqrt 3 s^* \geq r^*$.
Thus from the assumption $p \leq r^*$ we get $p \leq r$.

Translate $C$ by the vector $[-p, -q]$ and denote by $C'$ the image.  
Simultaneously, by this vector $[-p, -q]$, translate $H_C = c_1...c_6$ up to $H_{C'} = c'_1...c'_6$.  
We have $c'_i = (\cos \frac{2\pi}{6}i -p, \sin \frac{2\pi}{6}i -q)$, for $i=1, \dots, 6$ (see Fig. 1).

Let ${\overline{c}}_i = (\sqrt 3 \cos (\frac{2\pi}{6}i - \frac{\pi}{6}), \sqrt 3 \sin (\frac{2\pi}{6}i - \frac{\pi}{6}))$ for $i=1, \dots, 6$. 
Take the translation ${\overline{c}}'_i$ of ${\overline{c}}_i$ by $[-p, -q]$.
Thus ${\overline{c}}_i = (\sqrt 3 \cos (\frac{2\pi}{6}i - \frac{\pi}{6})-p, \sqrt 3 \sin (\frac{2\pi}{6}i - \frac{\pi}{6})-q)$.
Construct the star $S(H_{C'})$ over $H_{C'}$ being the union of the triangles $\overline c_1' \overline c_3', \overline c_5'$ and $\overline c_2' \overline c_4' \overline c_6'$. 

Take into account $\tau(D)$ and $\tau(H_{D} )= \tau(d_1) \dots \tau(d_6)$. 
We have $\tau(d_i) = (\sqrt 3 \cos (\frac{2\pi}{6}i - \frac{\pi}{6}), \sqrt 3 \sin (\frac{2\pi}{6}i -  \frac{\pi}{6}))$ for $i=1, \dots, 6$. 
Clearly $\tau(d_i) = \overline c_i$ for  $i=1, \dots, 6$.

{\ }

\vskip -0.9cm
\begin{figure}[htbp]
\hskip1.2cm \includegraphics[width=12.8cm,height=9.8cm]{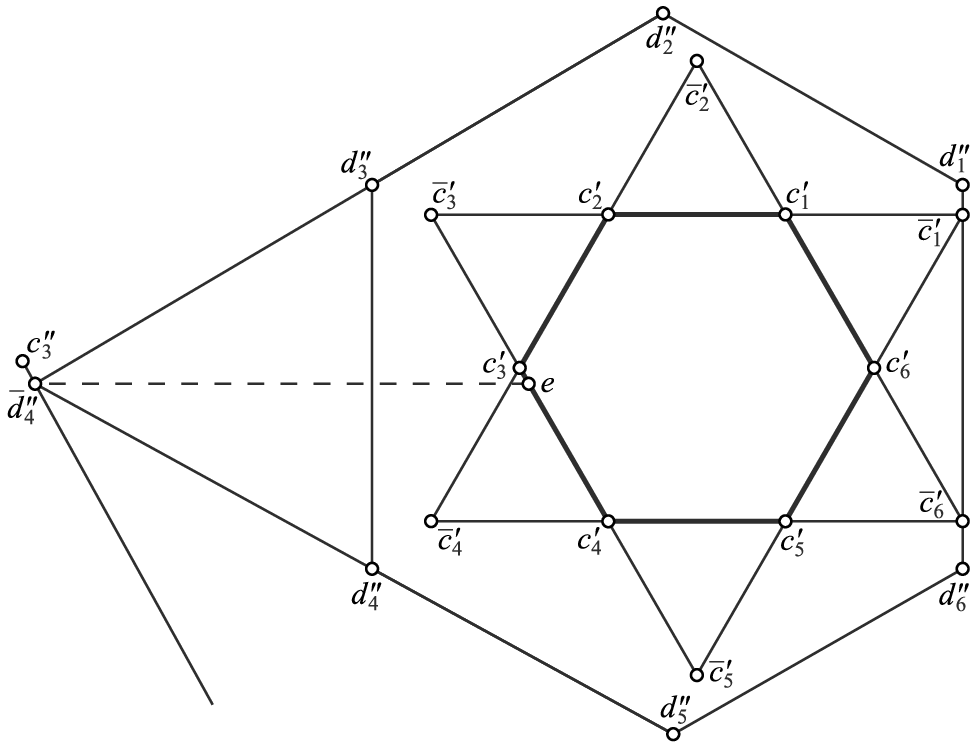} \\ 
\vskip-0.2cm
\caption{Illustration to the proof of Theorem} 
\end{figure}

\vskip-0.1cm

Translate $\tau(D)$ up to $D'$ and simultaneously $\tau(H_D)$ up to $H_{D'}$ by the vector $[-r, -s]$, which means that $H_{D'}$ has vertices $d'_i  = (\sqrt 3 \cos (\frac{2\pi}{6}i - \frac{\pi}{6}) -r, \sqrt 3 \sin (\frac{2\pi}{6}i -  \frac{\pi}{6}) -s)$ for $i=1, \dots, 6$. 

The bodies $C'$ and $D'$ have their centroids at $o = (0,0)$.  

Of course, $\delta_{\rm BM}^{\rm cen} (C', D') = \delta_{\rm BM}^{\rm cen} (C, D)$.

Since $0 \leq p \leq \frac{4}{21}$ and $0 \leq r \leq \frac{2}{7}$, then from $p \leq r$ we see that the number $\frac{3-2p}{3-2r}$ is at least~$1$.

Consider the homothetic image $H_{D''}$ of $H_{D'}$ with this ratio $\frac{3-2p}{3-2r}$.
Clearly, $d''_i = \frac{3-2p}{3-2r} d'_i$ for $i=1, \dots , 6$ are the vertices of $H_{D''}$ (again see Fig. 1). 
In particular, 
$d''_1 = (\frac{3-2p}{3-2r}(\frac{3}{2} -r),  \frac{3-2p}{3-2r}(\frac{1}{2}\sqrt 3 -s))$  
and
$d''_6 = (\frac{3-2p}{3-2r}(\frac{3}{2} -r),  \frac{3-2p}{3-2r}(-\frac{1}{2}\sqrt 3 -s))$.

We have $\frac{3}{2} - p = \frac{3-2p}{3-2r}(\frac{3}{2} -r)$ (in order to see this multiply both sides by $2(3-2r)$).
Thus the first coordinates of $d_1''$, $\overline c'_1$, $\overline c'_6$ and $d_6''$ are equal, i.e., these points are in one vertical line.
Moreover, by $p \leq r$ and $q \geq s$, we have $(3- 2r)(\sqrt 3 -2q) \leq (3-2p)(\sqrt 3 - 2s)$ which leads to $\frac{1}{2}\sqrt 3 - q \leq \frac{3-2p}{3-2r}(\frac{1}{2}\sqrt3 -s)$.
Consequently, ${\overline c'}_6{\overline c'}_1 \subset d''_6d''_1$.
Thus $H_{D''}$ contains ${\overline c'}_6$ and ${\overline c'}_1$.
We let the reader to show that $H_{D''}$ contains also the remaining ${\overline c'}_i$.  
Hence $H_{D''} \supset S(H_{C'})  \supset C'$.
Taking into account the inclusion $H_{D''} \subset D''$ we get $C' \subset D''$.

Prolonging $d_1''d_2''$, $d_3''d_4''$ and $d_5''d_6''$ we obtain a triangle.
Prolonging $d_2''d_3''$, $d_4''d_5''$ and $d_6''d_1''$ we also obtain a triangle.
The star being the union of these two triangles is denoted by $S(H_{D''})$ and called the star over $H_{D''}$
Denote its ``outer" vertices by ${\overline {d''}}_i$ for $i=1, \dots, 6$ so that $d''_{i-1} \in d''_{i-2}\overline {d''}_i$ (mod $6$), where $i=1 \dots 6$.
Since $d_3'' = (\frac{3-2p}{3-2r}(-\frac{3}{2}-r), \frac{3-2p}{3-2r}(\frac{1}{2}\sqrt 3 -s)$ and
$d_4'' = (\frac{3-2p}{3-2r}(-\frac{3}{2}-r), \frac{3-2p}{3-2r}(-\frac{1}{2}\sqrt 3 -s)$,
its length is $\frac{3-2p}{3-2r}\sqrt 3$ and 
the middle of the segment $d_3''d_4''$ is 
$(\frac{3-2p}{3-2r}(-\frac{3}{2}-r), -\frac{3-2p}{3-2r}s)$. 
This and the fact that the triangle $d''_3\overline {d''}_4 d''_4$ is regular imply that the distance between this middle and $\overline {d''}_4 $ is $\frac{3}{2}\cdot \frac{3-2p}{3-2r}$.
Hence $\overline {d''}_4 = (-\frac{(3-2p)(3+r)}{3-2r}, -\frac{(3-2p)s}{3-2r})$.

The equation of the straight line of the side $c'_3c'_4$ of $H_{C'}$ is $y= -\sqrt 3(x+1+p) -q$. 
Intersect it with the horizontal line $y= -\frac{(3-2p)s}{3-2r}$ passing through $\overline {d''}_4$. 
We obtain the point $e$ of intersection whose first coordinate is $\frac{\sqrt 3}{3} \cdot \frac{(3-2p)s}{3-2r} -1 -p -\frac{1}{3}\sqrt 3 q$. 
We see that $\overline {d''}_4$ and $e$ are on the same horizontal level.
By $f(p,q,r,s)$ denote the quotient of $-\frac{(3-2p)(3+r)}{3-2r}$ by $\frac{\sqrt 3}{3} \cdot \frac{(3-2p)s}{3-2r} -1 -p -\frac{1}{3}\sqrt 3 q$.
In other words, this is $|o\overline {d''}_4|/|oe|$.

Put $H_{C''} = f(p,q,r,s) \cdot H_{C'}$. 
Its vertices are $c_i'' = f(p,q,r,s) \cdot c'_i$ (in Figure 1 we see only $c_3''$).
Apply the intercept (Thales's) theorem for the three straight lines containing the segments $oc_3''$, $o\overline {d''}_4$ and $oc_3''$, and for the two parallel lines containing the segments $c''_3 c''_4$ and $c'_3 c'_4$.
We obtain $\overline {d''}_4 \in c''_3c''_4$
This and $c''_3c''_4 \subset H_{C''}$ imply $d{d''}_4 \in H_{C"}$.
We let the reader to show that also the remaining $\overline {d''}_i$ are in $H_{C''}$ and thus in $f(p,q,r,s) \cdot C'$.
Hence $S(H_{D''}) \subset f(p,q,r,s) \cdot C'$.
This and $D'' \subset S(H_{D''})$ imply that $D'' \subset f(p,q,r,s) \cdot C'$.
So by $C' \subset D''$ (established earlier in this proof) we conclude that for every $(p, q, r, s) \in T \times T^+$ (i.e., when $(p,q) \in T$ and $(r,s) \in T^+$) we get $\delta_{BM}^{\rm cen} (C', D') \leq f(p, q, r, s)$ and thus that

$$\delta_{BM}^{\rm cen} (C, D) \leq f(p, q, r, s).$$

Having in mind the description of $f(p,q,r,s)$ as a quotient is easy to evaluate it:

$$f(p,q,r,s) = \frac{\sqrt 3 (3-2p)(3+r)}{(3-2r)(\sqrt 3 + \sqrt 3 p + q) - (3-2p)s}.$$

It remains to show that $f(p, q, r, s) \leq \frac{69}{17}$ for every $(p, q, r, s) \in T \times T^+$.

From the formula for $f(p, q, r, s)$ we see that if $q$ decreases and $s$ increases while $p, r$ are constant, then $f(p,q,r,s)$ increases.
So the greatest value of $f(p, q, r, s)$ in $T \times T^+$ must be for some $(p,q)$ in the horizontal side of $T$ and for some $(r,s)$ on the horizontal side of $T^+$.

Consider the rectangle $Q = \{(p, r); \; 0 \leq p \leq \frac{4}{21}, 0 \leq r \leq \frac{2}{7} \}$.

Putting $q = 0$ and $s = 0$ into the above formula for $f(p,q,r,s)$ we get the function  

$$g(p, r) = \frac{\sqrt 3 (3-2p)(3+r)}{(3-2r)(\sqrt 3 + \sqrt 3 p)}.$$

Below, for this function, we apply the following method of finding the global maximum of a continuous function $f(x,y)$ in a polygon $Q \subset E^2$.
Namely, first we find the points being the solutions of the system of two equations when partial derivatives of our function $f(x,y)$ are $0$ in the interior of $Q$.  
Next we write the equations of the sides in the forms $y=u(x)$ or $x=v(y)$.
We find the critical points in the relative interiors of each side, where the derivative of the respective equation is $0$.
Finally, we check the values of $f(x,y)$ at the vertices of $Q$. 
The largest value at all the found points gives the maximum value of $f(x,y)$ in $Q$. 

We evaluate the partial derivatives $g'_p(p, r)$ and $g'_r(p, r)$ and easily check that the system of equations $g'_p(p, r) = 0$ and $g'_r(p, r) = 0$ has no solution in the interior of $Q$.

If we insert $0$ in place of $p$, we get a homographic function of the variable $r$.
The same if we insert $\frac{4}{21}$.
Also if we insert $0$ in place of $r$, we get a homographic function of the variable $p$.
The same if we insert $\frac{2}{7}$.
Since none of homographic functions attains a local extremum on an open interval,
we conclude that $g(p,r)$ considered on $Q$ does not attain a global extremum on the interiors on each side of $Q$. 

Finally, we check the values of $g(p, r)$ at the corners of the rectangle $Q$.
It appears that the greatest value is attained at the vertex 
$(0, \frac{2}{7})$ of $Q$.
Its value equals to $\frac{69}{17}$.

We conclude that the greatest value of $g(p, r)$ on $Q$ is  $\frac{69}{17}$.
Consequently, $\delta_{BM}^{\rm cen} (C, D) \leq \frac{69}{17}$. 
\end{proof}

\section{Final remarks}

On the turn of pages 259--260 of \cite{[G]}, Gr\"unbaum says that probably the least upper bound of 
$\delta_{\rm BM}^{\rm cen} (C, D)$ is $\frac{5}{2}$ (he uses a different notation).
 Without a proof, it is added that it can be shown that $\frac{5}{2}$ is attained for the parallelogram and triangle. 
A proof of this fact is given in \cite{[L5]}.

For comparison recall the conjecture from \cite{[L1]} that $\delta_{\rm BM} (C, D) \leq \frac{\cos^2 36}{\sin 18}$ which is  $\frac{2+ \sqrt 5}{2} \approx 2.118$.
This value is attained for the regular pentagon and the triangle (in their optimum positions the centroids do not coincide).

By the way, $\delta_{\rm BM}^{\rm cen} (P, T) \leq \frac{7- \sqrt 5}{2} \approx 2.382$. 
Namely, in Fig. 2 we see the regular pentagon $P$ with vertices $(\cos \frac{2\pi}{5}i, \sin \frac{2\pi}{5}i)$, where $i=1, \dots, 5$, and the inscribed triangle $T$ with vertex $(1,0)$, whose opposite side is in the straight line $x= -\frac{1}{2}$. 
We also see the homothetic image $\frac{7- \sqrt 5}{2} T$.
The centroids of $T$, $P$ and $\frac{7- \sqrt 5}{2} T$ coincide.
We suppose that $\delta_{\rm BM}^{\rm cen} (P, T) = \frac{7- \sqrt 5}{2}$ for the

\vskip0.2cm

\begin{figure}[htbp]
\hskip4.5cm \includegraphics[width=6.9cm,height=7.15cm]{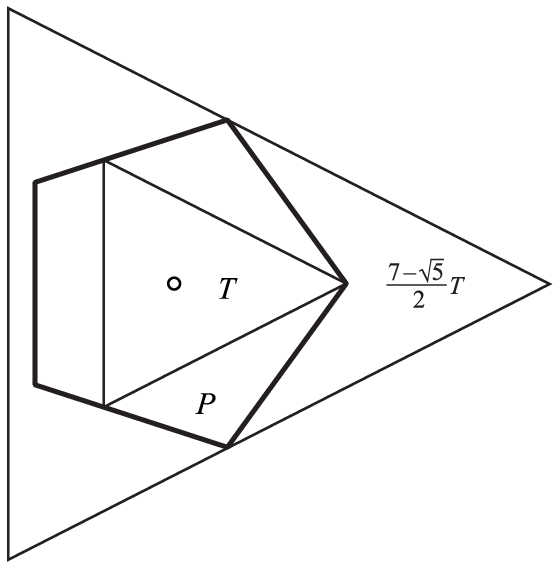} \\ 
\vskip-0.5cm
\caption{Illustration to the inequality $\delta_{\rm BM}^{\rm cen} (P, T) \leq \frac{7- \sqrt 5}{2}$} 
\end{figure}

\vskip0.15cm

\noindent 
regular pentagon $P$ and triangle $T$.
Is $\frac{7- \sqrt 5}{2}$ the supremum of $\delta_{\rm BM}^{\rm cen} (C,D)$ over all $C, D \in {\mathcal C}^2$ ?

At the end, let us show that $\delta_{\rm BM}^{\rm cen} (C,D) = \delta_{\rm BM}^{\rm cen} (D,C)$ for every $C, D \in {\mathcal C}^d$, as said in the Introduction. 
Take into account

$$\delta_{\rm BM}^{\rm cen} (D,C) = \inf_{a, h_\mu} \{\mu ; \, b(D) \subset C \subset h_\mu a(C) \ {\rm and} \ {\rm cen}(b(D)) =  {\rm cen}(C) \}.$$

Our considerations are not narrower if we assume that all the centroids in $\delta_{\rm BM}^{\rm cen} (C,D)$ and $\delta_{\rm BM}^{\rm cen} (D,C)$ are at one point.
Look at $a(C) \subset D \subset h_{\lambda} a(C)$ from the definition of $\delta_{\rm BM}^{\rm cen} (C,D)$. 
From $a(C) \subset D$ we get $h_{\lambda} a(C) \subset h_{\lambda} D$.
By this and by $D \subset h_{\lambda} a(C)$ we obtain $D \subset h_{\lambda} a(C) \subset h_{\lambda} D$.
Apply the affine transformation $b= h_{1/\lambda}a^{-1}$, which means $a^{-1} = h_\lambda b$, to these inclusions.
Hence $b(D) \subset C \subset h_{\lambda}b(D)$ is equivalent to $a(C) \subset D \subset h_{\lambda} a(C)$.
This leads to $\delta_{\rm BM}^{\rm cen} (D,C) = \delta_{\rm BM}^{\rm cen} (C,D)$.

\baselineskip 11 pt
\noindent
Marek Lassak

\noindent
University of Science and Technology

\noindent
al. Kaliskiego 7, 85-789 Bydgoszcz, Poland

\noindent
e-mail: lassak@pbs.edu.pl


\baselineskip 6 pt

\begin{thebibliography}{10}

\bibitem{[AS]}
G. Aubrun and S. J. Szarek, Alice and Bob Meet Banach, The interface of asymptotic geometric analysis and quantum information theory. Mathematical Surveys and Monographs, 223. American Mathematical Society, Providence, RI, 2017.

\bibitem{[Ba]} 
S. Banach, Th\'eorie des op\'erations lin\'eaires, Monogr. Mat. {\bf 1}. Warszawa (1932).  [English translation: Theory of linear operations. Translated from the French by F. Jellett. With comments by A. Pe{\l}czy\'nski and Cz. Bessaga. North-Holland Mathematical Library, 38. North-Holland Publishing Co., Amsterdam, 1987.] 

\bibitem{[Be]}  A. S. Besicovitch, {\it Measure of assymetry of convex curves}, J. London Math. Soc. {\bf 23} (1948), 237--240.  

\bibitem{[B+F]}
T. Bonnesen and W. Fenchel, Theorie der konvexen K\"orper, Springer, Berlin, 1934. English translation: Theory of convex bodies, BBC Associates, Moscow, Idaho USA, 1987

\bibitem{[G]} 
B. Gr\" unbaum, {\it Measures of symmetry of convex sets}, Convexity, Proc. Sympos. Pure Math., vol. 7, Amer. Math. Soc., Providence, R.I., 1963, pp. 233--270.

\bibitem {[L1]} 
M. Lassak, {\it Approximation of convex bodies by triangles}, Proc. Amer. Math. Soc. {\bf 115} (1992), 207--210.

\bibitem {[L3]} 
M. Lassak, {\it Banach-Mazur distance from the parallelogram to the affine-regular hexagon and other affine-regular even-gons}, Results Math. {\bf 76} (2021), 76--82.

\bibitem {[L4]}
M. Lassak, {\it Position of the centroid of a planar convex body}, (to appear in Aequationes Math.) temporarily see the manuscript arXiv:2202.01815v4

\bibitem {[L5]}
M. Lassak, {\it The centroid Banach-Mazur distance between the parallelogram and the triangle}, (to appear in J. Convex Analysis) temporarily see the manuscript arXiv:2210.16150v1.

\bibitem{[T-J]}
N. Tomczak-Jaegerman, Banach-Mazur Distances and Finite-dimensional Operator Ideals. Longman Scientific and Technical (Harlow and New York), 1989.

\bibitem{[T]}
G. Toth, Measures of Symmetry for Convex Sets and Stability, Universitext, Springer, Cham, 2015

\end{thebibliography}
\end{document}